\documentclass[12pt, reqno]{amsart}
\usepackage{amsmath, amsthm, amscd, amsfonts, amssymb, graphicx, xcolor}
\usepackage{amsmath}
\usepackage{cases}
\usepackage{amssymb}
\usepackage{mathrsfs}
\usepackage{euscript}
\usepackage{dsfont}
\usepackage{graphicx}
\usepackage{float}
\usepackage{graphicx}
\usepackage{amssymb}
\usepackage[]{cite}
\usepackage[bookmarksnumbered,colorlinks=true,plainpages=false,linkcolor=black,citecolor=black,urlcolor=black]{hyperref}

\allowdisplaybreaks[4]
\newtheorem{theorem}{Theorem}[section]

\newtheorem{lemma}{Lemma}[section]
\newtheorem{proposition}{Proposition}

\theoremstyle{definition}

\newtheorem{conjecture}[theorem]{Conjecture}

\theoremstyle{remark}
\newtheorem{remark}[theorem]{Remark}
\numberwithin{equation}{section}

\begin{document}
	\setcounter{page}{1}
	

	
	\title[]{Note on the largest prime factors of consecutive integers II}

	\author[]{Zhiyuan Yang}
	
	\address{School of Mathematics,  Shandong University, Jinan 250100, Shandong, China}
	\email{\tt zhiyuan.yang@mail.sdu.edu.cn}

	
	

	

	\begin{abstract}
		For an integer $n > 1$, let $P^+(n)$ be the largest prime factor of $n$. We show that a positive proportion of integers $n\leq x$ satisfy\begin{align*}
			P^+(n-1)<P^+(n)>P^+(n+1)<P^+(n+2)>P^+(n+3).
		\end{align*}
		The proof employs Selberg's upper-bound sieve and the multidimensional Selberg sieve developed in the Polymath8b work.
		\newline
		\newline
		\noindent \textit{Keywords.} largest prime factor; Selberg sieve.
		\newline
		\noindent 
	\end{abstract} \maketitle


	\section{Introduction}
	Let $P^+(n)$ denote the largest prime factor of $n$ with the convention that $P^+(1)=1$. In their correspondence in the 1930s, Erd\H{o}s and Tur\'an
	formulated the following conjecture; see S\'os~\cite{s}.
	\begin{conjecture}\label{conj:erdos-turan}
		(Erd\H{o}s-Turán) As $x\rightarrow\infty$, we have
		\begin{align*}
			\#\{n\leq x:P^+(n)<P^+(n+1)\}\sim \frac{1}{2}x.
		\end{align*}
	\end{conjecture}
	In 1978, Erd\H{o}s and Pomerance~\cite{ep} first proved that there exists a positive asymptotic density of integers $n$ such that $P^+(n)<P^+(n+1)$. In fact, they showed that\begin{align*}
		\#\{n\leq x:P^+(n)<P^+(n+1)\}>0.0099x.
	\end{align*}
	L\"u and Wang~\cite{lw} improved the constant to $0.2017$.
	The author's preprints~\cite{y1,y2} subsequently improved it to $0.280$
	and $0.299$, respectively. See~\cite{lw} for further historical background.
	
	One may also compare the largest prime factors of more than two consecutive
	integers. De Koninck and Doyon~\cite{kd} formulated the following
	generalization in 2011.
	\begin{conjecture}\label{conj:dekoninck-doyon}
		(De Koninck-Doyon) Fix an arbitrary integer $k\geq 2$ and let $n$ be a large number. Let $a_1,a_2,\cdots,a_k$ be any permutation of the numbers $0,1,\cdots,k-1$. Then we have\begin{align*}
			\#\{n\leq x:P^+(n+a_1)<P^+(n+a_2)<\cdots<P^+(n+a_k)\}\sim\frac{1}{k!}x.
		\end{align*}
	\end{conjecture}
	There has also been some progress for this type of comparison problem. 
	In 1978, Erd\H{o}s and Pomerance~\cite{ep} observed that
	the patterns
	\begin{align*}
		P^+(n-1)>P^+(n)<P^+(n+1),\ \text{or}\ P^+(n-1)<P^+(n)>P^+(n+1)
	\end{align*}
	and 
	\begin{align*}
		P^+(n-1)<P^+(n)<P^+(n+1)
	\end{align*}
	both occur infinitely often and conjectured that these two patterns occur for a
	positive density of $n$. In addition, they proved that there are infinitely many such integers. In 2001, Balog~\cite{ba} proved the following lower bounds
	\begin{align*}
		&\#\{n\leq x:P^+(n-1)<P^+(n)<P^+(n+1)\}\gg x^{1/2},\\&\#\{n\leq x:P^+(n-1)>P^+(n)>P^+(n+1)\}\gg x^{1/2}.
	\end{align*}
	Later, in 2018, Wang~\cite{w1} obtained positive-density results for two patterns by showing
	\begin{align*}
		&\#\{n\leq x:P^+(n-1)>P^+(n)<P^+(n+1)\}\gg x,\\&\#\{n\leq x:P^+(n-1)<P^+(n)>P^+(n+1)\}\gg x.
	\end{align*}
	For more difficult cases, Tao and Teräväinen~\cite{tt} succeeded in obtaining a positive
	density by showing
	\begin{align*}
		&\#\{n\leq x:P^+(n-1)<P^+(n)<P^+(n+1)>P^+(n+2)\}\gg x,\\&\#\{n\leq x:P^+(n-1)>P^+(n)>P^+(n+1)<P^+(n+2)\}\gg x.
	\end{align*}
	However, Tao and Teräväinen's proof is very complicated.
	In particular, using estimates on the numbers of prime $k$-tuples,
	Kim~\cite{k} proved the following pair of results
	\begin{align*}
		&\#\{n\leq x:P^+(n-1)<P^+(n)<P^+(n+1)\}\gg \frac{x}{(\log x)^{50}},\\&\#\{n\leq x:P^+(n-1)>P^+(n)>P^+(n+1)\}\gg \frac{x}{(\log x)^{50}}.
	\end{align*}
	Very recently, motivated by~\cite{dw}, the author~\cite{y3} gave an alternative proof of Tao and Teräväinen's results by showing
	\begin{align*}
		&\#\{n\leq x:P^+(n-1)<P^+(n)<P^+(n+1)>\max_iP^+(n+a_i)\}\gg x,\\&\#\{n\leq x:P^+(n-1)>P^+(n)<P^+(n+1)>\max_iP^+(n+a_i)\}\gg x,
	\end{align*}
	where the $a_i$ are finite in number, with $|a_i| = O(1)$. Moreover, using the results in~\cite{w2}, under the Elliott-Halberstam conjecture for friable integers, the author~\cite{y3} also obtains several comparison results for five numbers with positive density; one of them is
	\begin{align*}
		\#\{n\leq x:P^+(n-1)<P^+(n)>P^+(n+1)<P^+(n+2)>P^+(n+3)\}\gg x.
	\end{align*}
	In this paper, we reprove this result unconditionally.
	\begin{theorem}\label{thm:main}
		For sufficiently large $x$, we have
		\begin{align*}
			\#\{n\leq x:P^+(n-1)<P^+(n)>P^+(n+1)<P^+(n+2)>P^+(n+3)\}\gg x.
		\end{align*}
	\end{theorem}

	\section{Preliminaries}
	
	Let $\mu$, $\varphi$, and $\tau$ denote the M\"obius function,
	Euler's totient function, and the divisor function, respectively.
	We write $P^-(m)$ for the least prime factor of $m>1$, with $P^-(1)=\infty$,
	and $\Omega(m)$ for the number of prime factors counted with multiplicity.
	The notation $(u,v)$ and $[u,v]$ denotes the greatest common divisor and
	least common multiple, respectively, and $\log_x u:=\log u/\log x$.
	All linear forms have fixed integer coefficients and positive leading
	coefficients. Unless otherwise stated, the implied constants are independent of $\delta$.
	\begin{lemma}\label{lem:mean-values}
		Let $\varphi$ be Euler's totient function, and let $\psi$ be Dedekind's $\psi$-function, i.e.,
		\begin{align*}
			\varphi(n)=n\prod_{p|n}\left(1-\frac{1}{p}\right),\quad \psi(n)=n\prod_{p|n}\left(1+\frac{1}{p}\right).
		\end{align*} We have
		\begin{align*}
			&\sum_{m\leq x}\left(\frac{\varphi(m)}{m}\right)^k=c_kx+O((\log x)^k),\\&\sum_{m\leq x}\left(\frac{\psi(m)}{m}\right)^k=c_k'x+O((\log x)^{k})
		\end{align*}
		where
		\begin{align*}
			c_k=\prod_p\left(1+\frac{(p-1)^k-p^k}{p^{k+1}}\right),\quad c_k'=\prod_p\left(1+\frac{(p+1)^k-p^k}{p^{k+1}}\right).
		\end{align*}
	\end{lemma}
	\begin{proof}
		See \cite[Theorems 3.2 and 3.5]{ss}.
	\end{proof}
	
	\begin{lemma}\label{lem:selberg-upper}(Selberg’s upper-bound sieve) Let $\xi\geq 2$ and $L>1$. Let $\kappa$ be the sieve dimension. Let $\mathscr{P}$ denote a set of primes. Let $\nu(n)$ denote the number of distinct prime factors of $n$. Let $z\geq 2$ and write $P(z):=\prod_{p\leq z,p\in\mathscr{P}}p$. There exists a sequence $\{\lambda_d\}_{d=1}^{\infty}$ of real numbers, vanishing for $d>\xi$ or $\mu(d)=0$, satisfying $\lambda_1=1$, $|\lambda_d|\leq 3^{\nu(d)}$, $\mu*1\leq \lambda*1$, and such that
		\begin{align*}
			\sum_{d|P(z)}\lambda_d\frac{\omega(d)}{d}\leq\frac{1}{G(\xi,z)}= \prod_{\substack{p\leq z\\p\in \mathscr{P}}}\left(1-\frac{\omega(p)}{p}\right)\left(\frac{1}{\sigma(s)}+O\left(\frac{Ls^{2\kappa+1}}{\log z}\right)\right),
		\end{align*}
		uniformly for all multiplicative functions $\omega$ satisfying
		\begin{align*}
			&(i)\ 0<\frac{\omega(p)}{p}\leq 1-\frac{1}{A_1}\ (p\in\mathscr{P}),\\&(ii)-L\leq\sum_{w\leq p<z}\frac{\omega(p)}{p}\log p-\kappa\log \frac{z}{w}\leq A_2\  (2\leq w\leq z),
		\end{align*}
		where $A_1\geq 1$, $A_2\geq 1$ are some suitable constants, $s=\log \xi/\log z$ and $\sigma(s)$ is determined by the differential-difference equations
		\begin{equation}
			\begin{aligned}
				&\left\{ \begin{aligned}
					&s^{-\kappa}\sigma(s)=2^{-\kappa}e^{-\gamma\kappa}\Gamma(\kappa+1)^{-1},\quad&& 0\leq s\leq 2,\\&s(\sigma(s))'=\kappa\sigma(s)-\kappa\sigma(s-2),\quad&&s> 2.
				\end{aligned}\right.\nonumber
			\end{aligned}
		\end{equation}
	\end{lemma}
	\begin{proof}
		This is Selberg’s upper-bound sieve. See \cite[Lemma 6.1]{hr}.
	\end{proof}
	
	A $k$-tuple of linear forms $\{a_1n+a_1',\cdots,a_kn+a_k'\}$ is said to be $admissible$ if for any prime $p$ there is $n_p\in\mathbb{Z}$ such that $p\nmid \prod_{i=1}^k(a_in_p+a_i')$. We consider the pairwise nonproportional tuples with 
	\begin{align*}
		\prod_i a_i\neq 0,\quad \prod_{i<j}(a_ia_j'-a_ja_i')\neq 0.
	\end{align*}
	
	We require the following upper bound from Selberg's sieve.
	\begin{proposition}\label{prop:upper-bound}
		Let $\{a_1n+a_1',\cdots , a_kn+a_k'\}$ be a tuple of pairwise nonproportional linear forms.
		There exist $\varepsilon_0, \delta_0>0$ depending only on the tuple such that for any $0<\delta<\delta_0$, there are $\ll \delta^k x$ integers $x<n\leq 2x$ for which every form admits a factorization
		\begin{align*}
			a_in+a_i'=m_{i,1}m_{i,2},\qquad
			m_{i,1}<x^\delta,\qquad P^-(m_{i,2})>x^{\varepsilon_0}.
		\end{align*}
		The constant implied by $\ll$ depends on $\varepsilon_0$, $k$, $a_i$ and $a_i'$, but is independent of $\delta$. 
	\end{proposition}
	\begin{proof}
		For fixed small factors $m_{i,1}$, one applies Selberg's upper-bound
		sieve to the remaining factors. Selberg’s upper-bound sieve does not change the order of the estimate, so we may proceed by direct computation. We may also use Theorem 5.7 in \cite{hr} to simplify the proof and Lemma~\ref{lem:mean-values} to determine the coefficient.
	\end{proof}

	\section{Preparation for multidimensional Selberg sieves}
	
	In this section, our main aim is to prove the proposition below.
	\begin{proposition}\label{prop:lower-bound}
		Let $k\geq54$, and let $\{a_1n+a_1',\cdots , a_kn+a_k'\}$ be an admissible tuple of pairwise nonproportional linear forms. There exist $\varepsilon_0, \delta_0>0$ depending only on the tuple such that for any $0<\delta<\delta_0$, there are $\gg \delta^k x$ integers $x<n\leq 2x$ for which every form admits a factorization
		\begin{align*}
			a_in+a_i'=m_{i,1}m_{i,2},\qquad
			m_{i,1}<x^\delta,\qquad P^-(m_{i,2})>x^{\varepsilon_0},
		\end{align*}
		and at least two forms satisfy $P^+(a_in+a_i')>x^{1-\delta}$.
		The constant implied by $\gg$ depends on $\varepsilon_0$, $k$, $a_i$ and $a_i'$, but is independent of $\delta$. 
		
	\end{proposition}
	
	We introduce the notation used throughout this section.
	Let $w:=\log\log x$ and $W:=\prod_{p\leq w}p$.
	Let $e_i<x^{\delta }$ satisfy 
	\begin{align*}
		\mu(e_1\cdots e_k)\neq 0,\quad\min_i P^-(e_i)>\max\{a_i,a_i',a_ia_j'-a_ja_i':1\leq i\leq k,1\leq j\leq k\}.
	\end{align*}
	Write $W'=W(\prod e_i)^2/(W,\prod e_i)$.
	Let $b\mkern-7mu\pmod{W'}$ satisfy
	\begin{align*}
		a_ib+a_i'\equiv 0\mkern-15mu\pmod{e_i},\quad \left(\prod\frac{a_ib+a_i'}{e_i}, W'\right)=1.
	\end{align*}
	The Chinese remainder theorem allows us to choose $b$ locally. We only need to consider the case where $p$ is large. When $(p,e_1\cdots e_k)=1$, it suffices to exclude all $-a_i'\overline{a_i}\pmod{p}$, so there are $p-k$
	permitted residues. If $p|e_i$, $e_i|(a_ib+a_i')$ implies $b\equiv -a_i'\overline{a_i}\pmod{p}$, and there are $p-1$
	permitted residues to avoid one residue class so that $(a_ib+a_i')/e_i$ is divisible by $p$.
	Thus the number of choices of $b$ is 
	\begin{align*}
		\asymp \prod_{\substack{O(1)\ll p\leq w\\(p,\prod e_i)=1}}(p-k)\prod_{p|\prod e_i}(p-1)\gg \frac{W'}{e_1\cdots e_k}\left(\frac{\varphi(W')}{W'}\right)^k.\tag{3.1}\label{eq:residue-count}
	\end{align*}
	Note that $n\equiv b\mkern-7mu\pmod{W'}$ implies for $i\neq j$,
	\begin{align*}
		\left(\frac{a_in+a_i'}{e_i},\frac{a_jn+a_j'}{e_j}\right)=1,\quad \left(\frac{a_in+a_i'}{e_i},W'\right)=1.\tag{3.2}\label{eq:coprimality}
	\end{align*}
	Define 
	\begin{align*}
		\theta_{i}(a_in+a_i'):=\log(a_in+a_i')1_{(a_in+a_i')/e_i\ \text{is prime}}.
	\end{align*}
	For each compactly supported $F:[0,+\infty)\rightarrow \mathbb{R}$, let
	\begin{align*}
		S(F):=\sup\{x\geq 0:F(x)\neq 0\},\qquad \lambda_F(m):=\sum_{d\mid m}\mu(d)F(\log_x d).
	\end{align*}
	We set $S(0)=0$.
	
	\begin{lemma}\label{lem:convolution}
		Let $N$ be a natural number
		coprime to $W'$ with $\log N = O(\log^{O(1)}x)$. Let $F_1,\cdots , F_k, G_1,\cdots , G_k : [0, +\infty) \rightarrow \mathbb{R}$ be
		fixed smooth compactly supported functions. Define
		\begin{align*}
			\mathscr{C}(W',N):=\sum_{\substack{d_1,\cdots,d_k,d_1',\cdots,d_k'\\ [d_1,d_1'],\cdots,[d_k,d_k'],W', N\ \text{coprime}}}\prod_{j=1}^k\frac{\mu(d_j)\mu(d_j')F_j(\log_x d_j)G_j(\log_x d_j')}{[d_j,d_j']}
		\end{align*}
		Then we have
		\begin{align*}
			\mathscr{C}(W',N)=(c+o(1))(\log x)^{-k}\left(\frac{W'N}{\varphi(W'N)}\right)^k,
		\end{align*}
		where
		\begin{align*}
			c:=\prod_{j=1}^k\int_0^{\infty}F_j'(t_j)G_j'(t_j)\mathrm{d}t_j.
		\end{align*}
		The same claim holds if the denominators $[d_j,d_j']$ are replaced by $\varphi([d_j,d_j'])$.
	\end{lemma}
	\begin{proof}
		This follows immediately from \cite[Lemma 30]{po}. In the subsequent proof, we seem to also need
		\begin{align*}
			\sum_{\substack{d_1,\cdots,d_k,d_1',\cdots,d_k'\\ [d_1,d_1'],\cdots,[d_k,d_k'],W', N\ \text{coprime}}}\prod_{j=1}^k&\frac{\mu(d_j)\mu(d_j')F_j(\log_R d_j)G_j(\log_R d_j')}{[d_j,d_j']}\\&=(c+o(1))(\log R)^{-k}\left(\frac{W'N}{\varphi(W'N)}\right)^k,
		\end{align*}
		where $(\log x)^{-1/2}\ll\log R/\log x\ll 1$.
	\end{proof}

	\begin{lemma}\label{lem:unweighted}
		For each fixed $1 \leq i \leq k$, let $F_i, G_i : [0, +\infty) \rightarrow \mathbb{R}$ be fixed smooth compactly
		supported functions. Assume
		\begin{align*}
			\sum_{i=1}^k(S(F_i)+S(G_i))<1-2k\delta.
		\end{align*}
		Then we have
		\begin{align*}
			\sum_{\substack{x<n\leq 2x\\n\equiv b\mkern-15mu\pmod{W'}}}\prod_{i=1}^k\lambda_{F_i}((a_in+a_i')/e_i)\lambda_{G_i}((a_in+a_i')/e_i)=(c+o(1))\frac{x}{W'(\log x)^{k}}\left(\frac{W'}{\varphi(W')}\right)^k,
		\end{align*}
		where
		\begin{align*}
			c:=\prod_{j=1}^k\int_0^{\infty}F_j'(t_j)G_j'(t_j)\mathrm{d}t_j.
		\end{align*}
	\end{lemma}
	\begin{proof}
		The proof is similar to that of \cite[Theorem 20 (i)]{po}.
		The left-hand side may be expanded as
		\begin{align*}
			\sum_{\substack{d_1,\cdots,d_k,d_1',\cdots,d_k'\\{}[d_1,d_1'],\cdots,[d_k,d_k'],W'\ \text{coprime}}}\Bigg(\prod_{i=1}^k\mu(d_i)\mu(d_i')F_i(\log_x d_i)G_i(\log_x d_i')\Bigg)\sum_{\substack{x<n\leq 2x\\n\equiv b\mkern-15mu\pmod{W'}\\a_in+a_i'\equiv 0\mkern-15mu\pmod{[d_i,d_i']},\forall i}}1.
		\end{align*}
		For the inner summation, we have
		\begin{align*}
			\sum_{\substack{x<n\leq 2x\\n\equiv b\mkern-15mu\pmod{W'}\\a_in+a_i'\equiv 0\mkern-15mu\pmod{[d_i,d_i']},\forall i}}1=\frac{x}{W'[d_1,d_1']\cdots [d_k,d_k']}+O(1).
		\end{align*}
		The $O(1)$ error is acceptable, and the main term follows from Lemma~\ref{lem:convolution}.
	\end{proof}

	\begin{lemma}\label{lem:prime-weighted}
		Let $1\leq i_0\leq k$ be fixed, and for each fixed $1 \leq i \leq k$ distinct from $i_0$, let $F_i, G_i : [0, +\infty) \rightarrow \mathbb{R}$ be fixed smooth compactly
		supported functions. Assume
		\begin{align*}
			\sum_{1\leq i\leq k,i\neq i_0}(S(F_i)+S(G_i))<1/2-2k\delta.
		\end{align*}
		Then we have
		\begin{align*}
			\sum_{\substack{x<n\leq 2x\\n\equiv b\mkern-15mu\pmod{W'}}}\theta_{i_0}(a_{i_0}n+a_{i_0}')&\prod_{1\leq i\leq k,i\neq i_0}\lambda_{F_i}((a_in+a_i')/e_i)\lambda_{G_i}((a_in+a_i')/e_i)\\&=(c+o(1))(1+O(\delta))\frac{x}{W'(\log x)^{k-1}}\left(\frac{W'}{\varphi(W')}\right)^{k}+\mathcal{R}_{i_0},
		\end{align*}
		where the error $\mathcal{R}_{i_0}$ satisfies $\sum_{e_i,\forall i}\sum_{b}\mathcal{R}_{i_0}\ll x(\log x)^{-A}$ for any $A>0$,
		\begin{align*}
			c:=\prod_{1\leq i\leq k,i\neq i_0}\int_0^{\infty}F_{i}'(t_i)G_i'(t_i)\mathrm{d}t_i.
		\end{align*}
	\end{lemma}
	\begin{proof}
		The proof is similar to that of \cite[Theorem 19 (i)]{po} with $\theta=1/2$.
		Without loss of generality, we assume that $i_0 = k$. As in the proof of Lemma~\ref{lem:unweighted}, we may use the Chinese remainder theorem to concatenate the congruence conditions on $n$ into a single primitive congruence condition
		\begin{align*}
			a_kn+a_k'=e_kp\equiv b'\mkern-15mu\pmod{a_kq},
		\end{align*}
		where $q=W'[d_1,d_1']\cdots [d_{k-1},d_{k-1}']$.
		Then the left-hand side may be expanded as
		\begin{align*}
			&\sum_{\substack{d_1,\cdots,d_{k-1},d_1',\cdots,d_{k-1}'\\{}[d_1,d_1'],\cdots,[d_{k-1},d_{k-1}'],W'\ \text{coprime}}}\Bigg(\prod_{i=1}^{k-1}\mu(d_i)\mu(d_i')F_i(\log_x d_i)G_i(\log_x d_i')\Bigg)\\&\times\sum_{\substack{(a_kx+a_k')/e_k<p\leq (2a_kx+a_k')/e_k\\e_{k}p\equiv b'\mkern-15mu\pmod{a_kq}}}\log(e_kp).
		\end{align*}
		A standard argument gives
		\begin{align*}
			&\sum_{\substack{(a_kx+a_k')/e_k<p\leq (2a_kx+a_k')/e_k\\e_{k}p\equiv b'\mkern-15mu\pmod{a_kq}}}\log (e_kp)=\frac{1}{\varphi(a_kq/e_k)}\sum_{\substack{(a_kx+a_k')/e_k<p\leq (2a_kx+a_k')/e_k}}\log (e_kp)\\&\quad\quad+\Bigg(\sum_{\substack{(a_kx+a_k')/e_k<p\leq (2a_kx+a_k')/e_k\\e_{k}p\equiv b'\mkern-15mu\pmod{a_kq}}}\log (e_kp)-\frac{1}{\varphi(a_kq/e_k)}\sum_{\substack{(a_kx+a_k')/e_k<p\leq (2a_kx+a_k')/e_k}}\log (e_kp)\Bigg).
		\end{align*}
		We call the first term the main term and the second term the error term.
		
		We cannot directly control the magnitude of the error term. However, we only need to sum over all residue classes modulo each $p|e_i$, and then the congruence for the residue class modulo $p$ becomes $a_{i}n + a_i' \equiv 0\mkern-7mu\pmod{e_i}$ and $((a_{i}n + a_i')/e_i,e_i)=1$. (The residue classes contributed by all $p<w$ are $\ll \log^{O(1)} x$, hence negligible for the final result.) We then remove this coprimality condition via M\"obius inversion, and the error term becomes
		\begin{align*}
			&\ll (\log x)\max_{b''}\sum_{\substack{q\\e_1\cdots e_{k-1}|q}}\tau^{O(1)}(q)\Bigg|\sum_{\substack{p\in I_k\\p\equiv b''\mkern-15mu\pmod{q}}}1-\frac{1}{\varphi(q)}\sum_{\substack{p\in I_k}}1\Bigg|\\&\ \ +\max_{b''}\sum_{\substack{q\\e_1\cdots e_{k-1}|q}}\tau^{O(1)}(q)\Bigg|\sum_{\substack{p\in I_k\\p\equiv b''\mkern-15mu\pmod{q}}}\log p-\frac{1}{\varphi(q)}\sum_{\substack{p\in I_k}}\log p\Bigg|,
		\end{align*}
		where $I_k=((a_kx+a_k')/e_k, (2a_kx+a_k')/e_k]$.
		We then sum over all $e_i$, and the total error term becomes $O(x(\log x)^{-A})$ for any $A>0$ by the Bombieri-Vinogradov theorem. The above procedure corresponds to $\sum_{e_i,\forall i}\sum_b \mathcal{R}_{i_0}\ll x(\log x)^{-A}$ in our lemma.
		
		Noting that $\varphi(a_kq/e_k)=\varphi(q)a_k/e_k=\varphi(W')\varphi([d_1,d_1']\cdots[d_{k-1},d_{k-1}'])a_k/e_k$, the final main term becomes
		\begin{align*}
			&\frac{e_k}{\varphi(W')a_k}\sum_{\substack{d_1,\cdots,d_{k-1},d_1',\cdots,d_{k-1}'\\{}[d_1,d_1'],\cdots,[d_{k-1},d_{k-1}'],W'\ \text{coprime}}}\prod_{i=1}^{k-1}\frac{\mu(d_i)\mu(d_i')F_i(\log_x d_i)G_i(\log_x d_i')}{\varphi([d_i,d_i'])}\\&\quad\times\sum_{\substack{(a_kx+a_k')/e_k<p\leq (2a_kx+a_k')/e_k}}\log(e_kp).
		\end{align*}
		The inner summation is independent of all $d_i$ and $d_i'$. Since
		\begin{align*}
			\sum_{\substack{(a_kx+a_k')/e_k<p\leq (2a_kx+a_k')/e_k}}\log (e_kp)=(1+o(1))\frac{a_kx}{e_k}\left(\frac{\log e_k}{\log(x/e_k)}+1\right),
		\end{align*}
		by Lemma~\ref{lem:convolution},
		we complete the proof of Lemma~\ref{lem:prime-weighted}.
		
	\end{proof}

	\begin{lemma}\label{lem:moments}
		For each fixed $1\leq i\leq k$, let $F_i,G_i:[0,+\infty)\rightarrow \mathbb{R}$ be fixed smooth compactly supported functions. Let
		$r_1,\cdots , r_k \geq 1$ and $s_1,\cdots , s_k \geq 0$ be fixed natural numbers. Then we have\begin{align*}
			&\sum_{\substack{x<n\leq 2x\\n\equiv b\mkern-15mu\pmod{W'}}}\prod_{j=1}^k\left(|\lambda_{F_j}((a_jn+a_j')/e_j)|^{r_j}\tau((a_jn+a_j')/e_j)^{s_j}\right)\ll \frac{x}{W'(\log x)^k}\left(\frac{W'}{\varphi(W')}\right)^k.\tag{3.3}\label{eq:moments}
		\end{align*}
		Furthermore, if $1 \leq j_0 \leq k$ is fixed and $p_0$ is a prime with $p_0 \leq x^{1/(10k)}$, then we have \begin{align*}
			&\sum_{\substack{x<n\leq 2x\\n\equiv b\mkern-15mu\pmod{W'}}}\prod_{j=1}^k\left(|\lambda_{F_j}((a_jn+a_j')/e_j)|^{r_j}\tau((a_jn+a_j')/e_j)^{s_j}\right)1_{p_0|((a_{j_0}n+a_{j_0}')/e_{j_0})}\\&\ll \frac{\log_x p_0}{p_0}\frac{x}{W'(\log x)^k}\left(\frac{W'}{\varphi(W')}\right)^k.\tag{3.4}\label{eq:prime-divisor}
		\end{align*}
		Finally, we also have
		\begin{align*}
			&\sum_{\substack{x<n\leq 2x\\n\equiv b\mkern-15mu\pmod{W'}}}\prod_{j=1}^k\left(|\lambda_{F_j}((a_jn+a_j')/e_j)|^{r_j}\tau((a_jn+a_j')/e_j)^{s_j}\right)1_{P^-((a_{j_0}n+a_{j_0}')/e_{j_0})\leq x^\varepsilon}\\&\ll \varepsilon\frac{x}{W'(\log x)^k}\left(\frac{W'}{\varphi(W')}\right)^k,\tag{3.5}\label{eq:small-prime}
		\end{align*}
		where the constant implied by $\ll$ is independent of $\varepsilon$.
	\end{lemma}

	\begin{proof}
		We closely follow the proof of \cite[Proposition 14]{po}. Similar to \cite[(50)]{po}, we only need to prove that
		\begin{align*}
			&\sum_{\substack{x<n\leq 2x\\n\equiv b\mkern-15mu\pmod{W'}}}\prod_{j=1}^k\Bigg(\tau((a_jn+a_j')/e_j)^{O(1)}\prod_{p|(a_jn+a_j')/e_j}\min(\sigma\log_x p,1)\Bigg)\\&\ll \frac{x}{W'(\log x)^k}\left(\frac{W'}{\varphi(W')}\right)^k\sigma^{O(1)},
		\end{align*}
		for any $\sigma\geq 1$. 
		By factoring the prime divisors of $(a_{j}n + a_j')/e_j$, the left-hand side of the above inequality becomes
		\begin{align*}
			\ll\sum_{*}\Bigg(\prod_{j=1}^k2^{O(1+\Omega(d_j))}\prod_{p|d_j}\min(\sigma\log_x p,1)\Bigg)\sum_{**}1
		\end{align*}
		where the outer sum $\sum_{*}$ is over $d_1,\cdots , d_k \leq x^{1/(10k)}$ with $d_1,\cdots , d_k, W'$ coprime, and the
		inner sum $\sum_{**}$ is over $x < n \leq 2x$ with $n \equiv b \mkern-7mu\pmod{W'}$ and $a_jn + a_j' \equiv0\mkern-7mu\pmod{d_j}$ for each $j$, with
		$P^-((a_jn+a_j')/(e_jd_j)) \geq R$ for each $j$. (Denote $d:=d_1\cdots d_k$ and $R:=x^{1/(10k(1+\Omega(d)))}$.)
		For auxiliary sieve divisors $h_i,h_i'$ with $[h_1,h_1'],\ldots,[h_k,h_k'],dW'$ pairwise coprime, \eqref{eq:coprimality} gives (Here, the $h_i$ correspond to the $e_i$ in~\cite{po}.)
		\begin{align*}
			\sum_{\substack{x<n\leq 2x\\a_in+a_i'\equiv 0\mkern-15mu\pmod{d_i[h_i,h_i']},\forall i\\n\equiv b\mkern-15mu\pmod{W'}}}1=\frac{x}{dW'[h_1,h_1']\cdots[h_k,h_k']}+O(1),
		\end{align*}
		we need only follow the remainder of the proof in \cite{po}.
		
		Replacing $d_{j_0}$ by $[d_{j_0}, p_0] < x^{1/(5k)}$, the same argument gives \eqref{eq:prime-divisor} from \eqref{eq:moments}.
		Finally, \eqref{eq:small-prime} follows immediately from \eqref{eq:moments} when $\varepsilon > 1/(10k)$, and from \eqref{eq:prime-divisor} and Mertens’
		theorem when $\varepsilon \leq 1/(10k)$.
	\end{proof}
	
	\subsection{Proof of Proposition~\ref{prop:lower-bound}}
	Define
	\begin{align*}
		\nu(n)=\Bigg(\sum_{j=1}^Jc_j\lambda_{F_{j,1}}((a_1n+a_1')/e_1)\cdots \lambda_{F_{j,k}}((a_kn+a_k')/e_k)\Bigg)^2,
	\end{align*}
	where $c_j$ are real numbers to be determined, and $F_{j,i}$ are compactly supported functions to be determined satisfying
	\begin{align*}
		\max_{1\leq j\leq J}\sum_{i=1}^kS(F_{j,i})<\frac{1}{2}\left(\frac{1}{2}-2k\delta\right).
	\end{align*}
	Expanding the square and applying Lemma~\ref{lem:prime-weighted}, we obtain
	\begin{align*}
		\sum_{\substack{x<n\leq 2x\\n\equiv b\mkern-15mu\pmod{W'}}}\nu(n)\theta_{i_0}(a_{i_0}n+a_{i_0}')=&(\beta_{i_0}+o(1))(1+O(\delta))\frac{x}{W'(\log x)^{k-1}}\left(\frac{W'}{\varphi(W')}\right)^k\\&+\mathcal{R}_{i_0},
	\end{align*}
	where $\mathcal{R}_{i_0}$ satisfies $\sum_{e_i,\forall i}\sum_b\mathcal{R}_{i_0}\ll x(\log x)^{-100}$,
	\begin{align*}
		\beta_{i_0}:=\sum_{j=1}^{J}\sum_{j'=1}^{J}c_jc_{j'}F_{j,i_0}(0)F_{j',i_0}(0)\prod_{1\leq i\leq k,i\neq i_0}\int_0^{\infty}F_{j,i}'(t_i)F_{j',i}'(t_i)\mathrm{d}t_i
	\end{align*}
	Lemma~\ref{lem:unweighted} also gives
	\begin{align*}
		\sum_{\substack{x<n\leq 2x\\n\equiv b\mkern-15mu\pmod{W'}}}\nu(n)=(\alpha+o(1))\frac{x}{W'(\log x)^k}\left(\frac{W'}{\varphi(W')}\right)^k,
	\end{align*}
	where
	\begin{align*}
		\alpha:=\sum_{j=1}^J\sum_{j'=1}^Jc_jc_{j'}\prod_{i=1}^k\int_0^\infty F_{j,i}'(t_i)F_{j',i}'(t_i)\mathrm{d}t_i.
	\end{align*}
	Let $a:=\max |a_i|$. If $\beta_1+\cdots+\beta_k>\alpha$, for sufficiently small $\delta_0$ and each fixed $\delta<\delta_0$, we have
	\begin{align*}
		\sum_{\substack{x<n\leq 2x\\n\equiv b\mkern-15mu\pmod{W'}}}\nu(n)\Bigg(\sum_{i=1}^k\theta_i(a_in+a_i')-\log(3ax)\Bigg)\geq& \left(2\beta_0+o(1)\right)\frac{x}{W'(\log x)^{k-1}}\left(\frac{W'}{\varphi(W')}\right)^k\\&+\sum_{i=1}^k\mathcal{R}_i,
	\end{align*}
	where $\beta_0=(\beta_1+\cdots+\beta_k-\alpha)/4$. (In fact, from the proof of \cite[Theorem 22]{po} and by \cite[Theorem 23 (vii)]{po}, when $k\geq 54$, we may find suitable $F_{j,i}$ such that $\beta_1+\cdots+\beta_k>\alpha$.)
	By \eqref{eq:small-prime}, we also have
	\begin{align*}
		\sum_{j_0=1}^k\sum_{\substack{x<n\leq 2x\\n\equiv b\mkern-15mu\pmod{W'}}}\nu(n)(\log x)1_{P^-((a_{j_0}n+a_{j_0}')/e_{j_0})\leq x^\varepsilon}\ll \varepsilon \frac{x}{W'(\log x)^{k-1}}\left(\frac{W'}{\varphi(W')}\right)^k,
	\end{align*}
	so for sufficiently small $\varepsilon_0$ (independent of $\delta$),
	\begin{align*}
		&\sum_{\substack{x<n\leq 2x\\n\equiv b\mkern-15mu\pmod{W'}}}\nu(n)\Bigg(\sum_{i=1}^k\theta_i(a_in+a_i')-\log(3ax)\Bigg)1_{\min_j P^-((a_{j}n+a_{j}')/e_{j})>x^{\varepsilon_0}}\\&\geq \left(\beta_0+o(1)\right)\frac{x}{W'(\log x)^{k-1}}\left(\frac{W'}{\varphi(W')}\right)^k+\sum_{i=1}^k\mathcal{R}_i.
	\end{align*}
	When $\min_jP^-((a_{j}n+a_{j}')/e_{j})>x^{\varepsilon_0}$, we have $|\nu(n)|=O(1)$, so the above result becomes
	\begin{align*}
		&\#\Big\{x<n\leq 2x:n\equiv b\mkern-15mu\pmod{W'},\sum_{i=1}^k\theta_i(a_in+a_i')>\log(3ax),\\&\quad\quad\quad\min_j P^-((a_{j}n+a_{j}')/e_{j})>x^{\varepsilon_0}\Big\}\\&\geq \left(\beta_0'+o(1)\right)\frac{x}{W'(\log x)^{k}}\left(\frac{W'}{\varphi(W')}\right)^k+\frac{1}{C\log x}\sum_{i=1}^k\mathcal{R}_i,
	\end{align*}
	where $\beta_0'>0$, $C>0$ are some suitable constants.
	Summing over all $b$ and $e_i$, by \eqref{eq:residue-count}, we have
	\begin{align*}
		&\sum_{e_i<x^\delta,\forall i}\#\Big\{x<n\leq 2x:e_i| (a_in+a_i'),\sum_{i=1}^k\theta_i(a_in+a_i')>\log(3ax),\\&\quad\quad\quad\min_j P^-((a_{j}n+a_{j}')/e_{j})>x^{\varepsilon_0}\Big\}\\&\gg \frac{x}{(\log x)^k} \sum_{\substack{e_i<x^\delta,\forall i\\\mu(e_1\cdots e_k)\neq 0}}\frac{1}{e_1\cdots e_k} +O\left(\frac{x}{(\log x)^{100}}\right)\\&\gg \delta^kx.
	\end{align*}
	Choose $\delta_0<\varepsilon_0$. Each  $e_i$ then
	contains exactly the prime-power factors of $a_in+a_i'$ below
	$x^{\varepsilon_0}$, so an integer cannot be counted for two different
	cofactor tuples. We complete the proof.
	
	\begin{remark}\label{rem:dependencies}
		For clarity, we explain the dependencies among the parameters: $\delta_0$ and $\varepsilon_0$ and the constant implied by $\gg$ depend at most on $k$, $a_i$, $a_i'$, and $F_{j,i}$, and are independent of $\delta$.
	\end{remark}
	
	\section{Proof of Theorem~\ref{thm:main}}
	
	\begin{lemma}\label{lem:coefficients}
		For any $k$, there exists $\alpha_1<\cdots<\alpha_k$ satisfying
		\begin{align*}
			(\alpha_j-\alpha_i)|\alpha_i,\quad \text{for}\ i<j.
		\end{align*}
	\end{lemma}
	\begin{proof}
		A related coefficient construction appears in \cite[Lemma 6 and Remark 9]{k}.
		Start with $\{1\}$. Suppose the set $\{\alpha_1,\dots,\alpha_s\}$ has been constructed. Let $D$ be a common multiple of all $\alpha_i$ and all positive differences $\alpha_j-\alpha_i$. Replace the original set with
		$\{D, D+\alpha_1, \dots, D+\alpha_s\}$. The property is preserved, and the induction is complete.
	\end{proof}
	
	Fix $k=54$, and let $D_0$ be a common multiple of all $\alpha_j-\alpha_i$. Denote
	\begin{align*}
		a_i:=\alpha_i+TD_0,
	\end{align*}
	where $T$ is a sufficiently large parameter to be determined.
	For $i<j$, write
	\begin{align*}
		d_{ij}=a_j-a_i,\quad C_{ij}=\frac{2a_ia_j}{d_{ij}},\quad B_{ij,s}=\frac{2a_i}{d_{ij}}+s,\quad s\in\{-1,1,3\}.
	\end{align*}
	By Remark~\ref{rem:large-shift}, we may choose $T$ sufficiently large such that
	\begin{align*}
		\Delta_{ij,s,\ell}:=C_{ij}-a_{\ell}B_{ij,s}\neq 0\quad (i<j,\ s=-1,1,3,\ 1\leq \ell\leq k).
	\end{align*}
	Consider $k$ admissible linear forms
	\begin{align*}
		L_\ell(n)=a_{\ell}n+1.\tag{4.1}\label{eq:linear-forms}
	\end{align*}
	From Proposition~\ref{prop:lower-bound}, we know that there are $\gg \delta^k x$ values of $x<n \leq 2x$ such that at least two of the linear forms in \eqref{eq:linear-forms} have largest prime factor exceeding $x^{1-\delta}$, while each of the remaining linear forms can be written as $m_1m_2$ with $m_1< x^\delta$ and $P^-(m_2) > x^{\varepsilon_0}$. Hence there exists at least one pair $i'<j'$ such that there are $\geq v_0\delta^k x$ values of $x<n\leq 2x$ for which $L_{i'}(n)$ and $L_{j'}(n)$ have largest prime factor exceeding $x^{1-\delta}$.
	Let $a=\max |a_i|$. Let 
	\begin{align*}
		N:=\frac{2a_{i'}}{d_{i'j'}}L_{j'}(n),
	\end{align*}
	then 
	\begin{align*}
		N+2=\frac{2a_{j'}}{d_{i'j'}}L_{i'}(n).
	\end{align*}
	Now we have
	\begin{align*}
		&\sum_{\substack{N\leq 5a^2x\\P^+(N-1)<P^+(N)>P^+(N+1)<P^+(N+2)>P^+(N+3)}}1\\&\geq v_0\delta^kx-\sum_{s=-1,1,3}\sum_{\substack{x<n\leq 2x\\e_i<x^\delta,e_i|L_i(n),P^-(L_i(n)/e_i)>x^{\varepsilon_0},\forall i\\P^+(2a_{i'}L_{j'}(n)/d_{i'j'}+s)>x^{1-\delta}}}1
	\end{align*}
	Fix $s$. By Proposition~\ref{prop:upper-bound}, we also know that
	\begin{align*}
		\sum_{e_i<x^\delta,\forall i}\sum_{e<x^{\delta}}\sum_{\substack{x<n\leq 2x\\e_i\mid L_i(n),\ (L_i(n)/e_i,P(x^{\varepsilon_0}))=1,\forall i\\e\mid 2a_{i'}L_{j'}(n)/d_{i'j'}+s\\((2a_{i'}L_{j'}(n)/d_{i'j'}+s)/e,P(x^{\varepsilon_0}))=1}}1\ll \delta^{k+1}x.
	\end{align*}
	Thus for a large $C$ (independent of $\delta$), we obtain
	\begin{align*}
		\sum_{\substack{N\leq 5a^2x\\P^+(N-1)<P^+(N)>P^+(N+1)<P^+(N+2)>P^+(N+3)}}1\geq (v_0-C\delta)\delta^kx.
	\end{align*}
	Taking $\delta<v_0/C$ completes the proof.

	\begin{remark}\label{rem:nonproportional}
		The two linear forms $an+b$ and $cn+d$ are not proportional, which is equivalent to
		\begin{align*}
			ad-bc\neq0.
		\end{align*}
		Thus for $s=-1,1,3$, $C_{ij}-a_\ell B_{ij,s}\neq0$ implies 
		\begin{align*}
			\frac{2a_{i}}{d_{ij}}L_{j}(n)+s\quad \text{and}\quad L_{\ell}(n)
		\end{align*}
		are not proportional.
	\end{remark}

	\begin{remark}\label{rem:large-shift}
		Substituting $a_r = A + \alpha_r$, we obtain
		\begin{align*}
			C_{ij}-a_{\ell}B_{ij,s}=0&\Leftrightarrow 2a_ia_j-2a_ia_{\ell}-sa_{\ell}(a_j-a_i)=0\\&\Leftrightarrow 2(A+\alpha_i)(A+\alpha_j)-2(A+\alpha_i)(A+\alpha_{\ell})-s(A+\alpha_{\ell})(\alpha_j-\alpha_i)=0
		\end{align*}
		This is either a linear equation or a constant equation in $A$. If the coefficient of the linear term is zero, the constant term equals
		\begin{align*}
			-\frac{s(2-s)}{2}(\alpha_j-\alpha_i)^2,
		\end{align*}
		which is non-zero for $s = -1, 1, 3$. Therefore each such equation excludes at most one value of $A$, and there are only finitely many such equations altogether.
	\end{remark}

	\section*{Acknowledgements}
	The author thanks ChatGPT as assistants for writing up the manuscript, checking, and typesetting some routine details.
	All new ideas in the paper were
	conceived without the use of AI. The author take full responsibility for the correctness.

\end{document}